\documentclass[12pt,a4paper]{article}
 
\usepackage{amsmath,amssymb,amsfonts}
\usepackage{epsfig}

\usepackage[latin1]{inputenc}
\usepackage[T1]{fontenc}
\usepackage{ae}
\usepackage[english]{babel}
 
\usepackage{fancyhdr}
\usepackage{lastpage}
\usepackage{geometry}
\numberwithin{equation}{section}

\newcounter{item}

\renewcommand{\theitem}{\arabic{section}.\arabic{item}}

\newcommand{\cc}{\setcounter{equation}{0}}

\newenvironment{theo}[1]{
\setcounter{item}{\value{equation}} 
\addtocounter{equation}{1}
\refstepcounter{item}
\par\addvspace{\bigskipamount}
\indent {\bf \theitem.\hspace{1em}Theorem#1.} \sl }
{\par\addvspace{\bigskipamount}
}

\newenvironment{pf}{
\par\addvspace{-\smallskipamount}
\indent {\bf Proof.}$\,\ $ } { $\Box$ 
\par\addvspace{\bigskipamount}  
}

\newenvironment{cor}{
\setcounter{item}{\value{equation}} \addtocounter{equation}{1}
\refstepcounter{item}
\par\addvspace{\bigskipamount}
\indent {\bf \theitem.\hspace{1em}Corollary.\,\ }  \sl }
{\par\addvspace{\bigskipamount}
}

\newenvironment{rem}{
\setcounter{item}{\value{equation}}\addtocounter{equation}{1}
\refstepcounter{item}
\par\addvspace{\bigskipamount}
\indent {\bf \theitem.\hspace{1em}Remark.$\, \ $} }
{\par\addvspace{\bigskipamount}
}

\newenvironment{nonsec}{\bf 
\setcounter{item}{\value{equation}}\addtocounter{equation}{1}
\refstepcounter{item}
\par\addvspace{\bigskipamount} 
\indent \theitem.\hspace{1em}\ignorespaces }
{\unskip .\ \ \ }

\newcounter{minutes}
\divide\time by 60
\newcounter{hours}
\multiply\time by 60
\addtocounter{minutes}{-\time}

\begin{document}

\begin{center}
{\Large \bf On conformal moduli of polygonal quadrilaterals}
\end{center}
\medskip

\begin{center}
{\large \bf  V.N.~Dubinin and M.~Vuorinen}
\end{center}
\bigskip
\begin{center}
\texttt{File:~\jobname .tex, 2005-06-28 - 2007-01-02 
        printed: \number\year-\number\month-\number\day, 
        \thehours.\ifnum\theminutes<10{0}\fi\theminutes}
\end{center}

\medskip

{{\bf Abstract.}
The change of conformal moduli of polygonal quadrilaterals
under some geometric transformations is studied. We consider
the motion of one vertex when the other vertices remain fixed,
the rotation of sides, polarization, symmetrization, and averaging
transformation of the quadrilaterals. Some open problems are
formulated.}

{{\bf 2000 Mathematics Subject Classification.} Primary 30C85. 
Secondary 30C75.}

\bigskip

\section{Introduction}{}

\smallskip

The notion of the conformal modulus of a quadrilateral has found 
very important applications to several questions in Geometric Function
Theory and its applications (cf. \cite{Ah}, \cite{AB}, 
\cite{LV}, \cite{K} and the references
therein). A simply connected domain of hyperbolic type in the extended
complex plane with four distinct marked accessible boundary points
$a,b,c,d$ (vertices) is called a quadrilateral. The points define, in
this order, a positive orientation with respect to the domain.
Following \cite{H}, we denote the quadrilateral by $(Q;a,b,c,d)$ and call 
its sides the boundary arcs between neighbouring vertices, the
vertices themselves excluded, and denote these by
$(a,b),(b,c),(c,d)$ and $(d,a)\,,$ respectively. Let the function
$w=f(z)$ be a one-to-one conformal map of the domain $Q$ onto the
rectangle $0<u<1, 0< v<M (w=u+iv)$ with the vertices $a,b,c,d$
corresponding to $0,1,1+iM, iM \,.$ The number $M$ is called the (conformal)
modulus of the quadrilateral $(Q;a,b,c,d)$ and we will denote it
$M(Q;a,b,c,d) \,. $  Note that $M(Q;a,b,c,d) = 1/M(Q;b,c,d,a) \,. $ The
modulus of a quadrilateral is a conformal invariant and agrees with
the extremal length of the family of curves joining 
the sides $(a,b)$ and $(c,d)$ in the domain $Q \,$ \cite{Ah}. 
In Physics, the modulus
means, for example, the reciprocal electrical resistance (up to a constant
multiple) of $Q$ as a metallic plate or an electrical conductor with
electrodes $\overline{(b,c)}$ and $\overline{(d,a)}$ 
(with a constant potential there). We will
also need the definition of the modulus of a quadrilateral as the a capacity
of a condenser \cite{PS}. In what follows we write
$$ I(v,B) := \iint_B |\nabla v|^2 dx dy \, .$$  
It is well-known that  
\begin{equation}
\label{eq1.1}
M(Q;a,b,c,d)= \min I(v,Q),
\end{equation}
where the minimum is taken over all admissible functions $v,$ i.e. real-valued
functions, continuous in $\overline{Q}\, ,$ satisfying a Lipschitz condition
in a neighborhood of every finite point of $Q,$ equal to zero on the side
$(d,a)$ and equal to one on the side
$(b,c)\,.$ The function $u(z) = Re f(z)$ is called the potential function of
the quadrilateral $(Q;a,b,c,d)\,.$ In the case of smooth sides
$(a,b)$ and $(c,d)$ the potential function may be characterized as the
function continuous in $\overline{Q}\,,$ harmonic in $Q\,,$ equal to zero
on $(d,a),$ equal to one on $(b,c)$ and such that the normal derivative
$\partial u / \partial n = 0$ on the other sides of the quadrilateral. From
the Dirichlet principle it follows that
\begin{equation}
\label{eq1.2}
M(Q;a,b,c,d)=  I(u,Q)\,,
\end{equation}
and $M(Q;a,b,c,d)<  I(v,Q)\,,$ if $v  \not\equiv u$ (cf. \cite[p. 434]{H}).
The minimum in  (\ref{eq1.1}) is called the capacity of the condenser with
the plates $\overline{(d,a)}$ and  $\overline{(b,c)}$ and the field $Q \,.$
More general condensers are studied in \cite{D1}. Note that the next two 
monotonicity properties of the modulus of a quadrilateral follow easily
from  (\ref{eq1.1})
and (\ref{eq1.2}) (cf. \cite[p.40]{D1}, \cite[p.436]{H}, \cite[p.128]{B}).

\medskip

{\bf Property 1.} {\sl If for two quadrilaterals $(Q;a,b,c,d)$ and 
$(Q;a',b',c',d')$
the inclusions $(b,c) \subset (b',c')$ and  $(d,a) \subset (d',a')$ hold, then
$$ M(Q;a,b,c,d) \le M(Q;a',b',c',d')$$
with equality if and only if the quadrilaterals coincide. }

\medskip

{\bf Property 2.} {\sl  If for two quadrilaterals $(Q;a,b,c,d)$ and 
$(Q';a,b,c,d)$ the domain $Q'$ is an extension of $Q$ through 
the sides $(a,b)$ and $(c,d)$, then
$$ M(Q;a,b,c,d) \le M(Q';a,b,c,d)$$
with equality only in the case $Q=Q' \,.$ }

\medskip

A consequence of the second property is that extension of the domain $Q$
across the sides $(b,c)$ and $(d,a)$ decreases the modulus of the
quadrilateral $(Q;a,b,c,d)\,.$ In what follows, unless otherwise stated,
we will understand that the quadrilaterals are polygonal, i.e. quadrilaterals
with (linear) intervals as sides. In this case the vertices $a,b,c,d$
uniquely define a domain $Q\,$ however, it will be convenient to adhere to
the above notation $(Q;a,b,c,d)\,.$
In spite of the considerable interest
of the modulus of quadrilateral, and the apparent simplicity of 
polygonal quadrilateral,
the behavior of the conformal modulus under geometric transformations
has been investigated very little in the literature
(cf. \cite[p.428]{H}, \cite[p.115]{K}, \cite[p.188]{PS}, \cite[p.128]{B},
\cite{S}).
The purpose of this work is to partially fill this gap. The motivation of this
paper is also connected with the well-known 
Schwarz-Christoffel formula \cite{AQVV}.
Numerical implementation of this formula in \cite{HVV} enabled one
to explore the values of the modulus of a quadrilateral and led to a number
of conjectures which will be settled here. In particular, we study the 
variation of the modulus of a
quadrilateral when one of its vertices is moving, the decrease of the modulus
under polarization, symmetrization, and averaging transformation of the
quadrilaterals.

\bigskip

\cc
\section{The motion of one vertex}{}\label{sec2}

It is clear that every transformation of a quadrilateral may be
considered as a motion of its vertices. In this section we consider 
simplest motions of one of the vertices, when the other vertices remain
fixed.

\begin{theo}{}\label{th2.1}
Let $(Q;a,b,c,d)$ be a quadrilateral for which the interior angle at the
vertex $a$ is at most equal to $\pi$ and let the point $a' $ be a point
of the side $(a,b)$ such that the segment $(d,a') \subset Q \,.$
Then
\begin{equation}
\label{eq2.2}
M(Q;a,b,c,d) < M(Q';a',b,c,d) \,.
\end{equation}
If the quadrilateral $(Q;a,b,c,d)$ is convex and the point $a', a' \neq a,$
lies in the closed convex set, whose boundary consists of the side $(a,b)$ and
the linear extensions of the sides $(d,a)$ and $(b,c)$ over the points $a,b,$
respectively, up to the point of their intersection 
(which may be $\infty$), then the inequality (\ref{eq2.2}) is again valid.
\end{theo}



\bigskip

\begin{figure}[ht]
\centerline{\psfig{figure=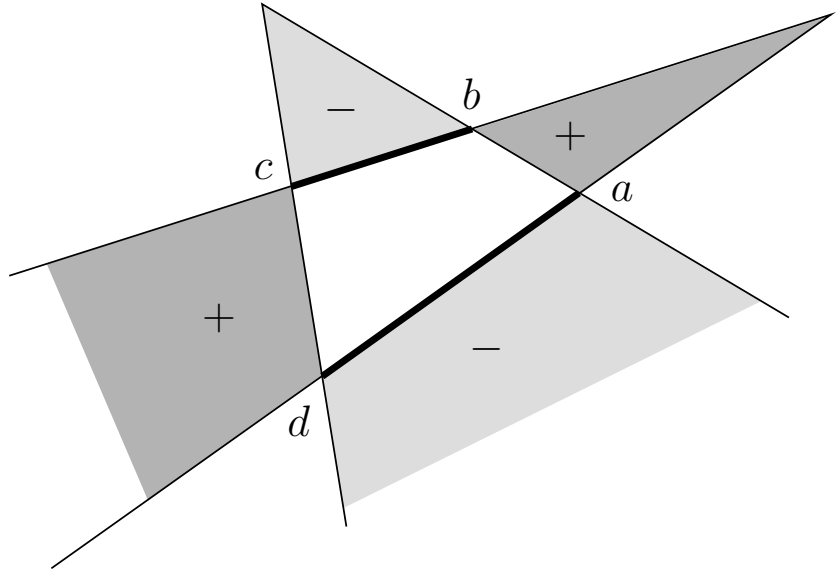,width=11cm}}
\caption{The modulus increases (decreases), when a vertex moves into a region
marked with + (-). }
\end{figure}

\begin{pf}
In the first case we consider the nonpolygonal quadrilateral $(Q;a',b,c,d)$ 
and the polygonal  quadrilateral $(Q';a',b,c,d) \,.$ Applying Properties 1
and 2 we have
$$ M(Q;a,b,c,d) < M(Q;a',b,c,d) < M(Q';a',b,c,d) \,.$$ 
In the second case we denote by $Q''$ the domain obtained by attaching
to the domain $Q$ the triangle $a a' b$ with its side $(a,b) \,.$ Then
$(Q'';a,b,c,d)$ is not polygonal. Again by  Properties 1 and 2
$$ M(Q;a,b,c,d) < M(Q'';a,b,c,d) < M(Q'';a',b,c,d) $$
$$ < M(Q';a',b,c,d) \,.$$ 
\end{pf}

Figure 1 displays the four convex sets, each bounded by a side of 
the quadrilateral and by the linear extensions of its two neighboring 
sides, and marked by $+$ or $- \,.$  In accordance
with Theorem \ref{th2.1} the modulus decreases (increases) 
if the vertex moves from its current location into
a region marked with $-$ ($+$ \, ).


\bigskip

\begin{figure}[ht]
\centerline{\psfig{figure=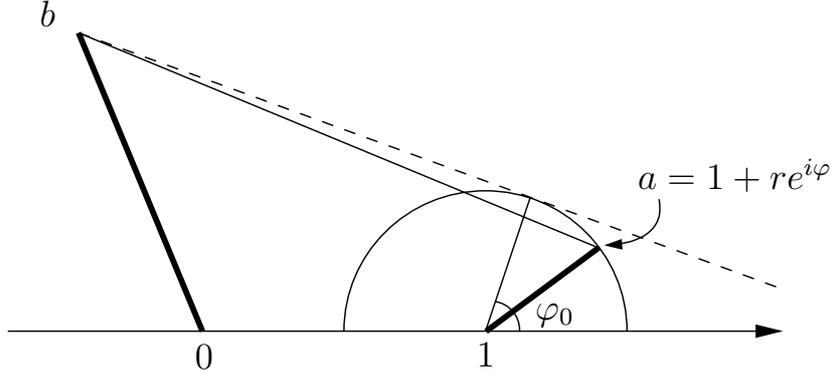,width=11cm}}
\caption{The modulus strictly increases on $[0, \varphi_0] .$}
\end{figure}

\begin{cor}
Let $0<r<1, Re b <1+r, Im b >0,  |b - 1|> r,$ and let $1+ r e^{i \varphi_0}$ be the point of 
tangency of the circle $|z-1|= r\, ,$ closest to $1+r,$  
with the line going through the point 
$b\,.$ Then the function $M(Q;1+r e^{i \varphi},b, 0, 1)$ strictly increases
on $[0, \varphi_0] \,.$
\end{cor}

\bigskip
\begin{figure}[ht]
\centerline{\psfig{figure=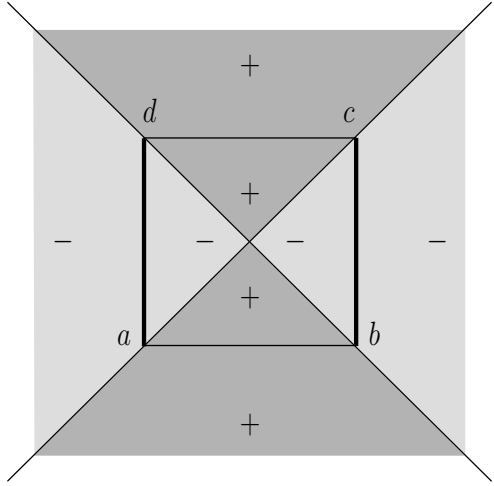,width=6.5cm,height=6.5cm}}
\caption{The regions of increase (+) and decrease (-) of modulus
in the case of a square.}
\end{figure}

In the particular case, when $Q$ is a square, Theorem \ref{th2.1} 
together with the well-known observation that the modulus of 
a quadrilateral, symmetric with respect to a diagonal, 
is equal to one \cite[p. 433]{H}, \cite[p. 218]{Her} give a complete
description of those sets, into which the displacement of one of the
vertices increases or decreases the numerical value of 
the conformal modulus of the square (Figure 3).

In the general case, the description of such sets seems to be difficult.
We only give the following result.

\begin{theo}{}\label{th2.3}
For fixed $M>0, \varphi \in (0, \pi/2)$ we have the formula
\begin{equation}
\label{eq2.4}
M(Q(\varphi,h);M,M+i,i, h e^{i \varphi }) = M + \frac{1}{2}(M \sin \varphi -
\cos \varphi)h + o(h)
\end{equation}
when $h \to 0 \,.$
\end{theo}

\begin{pf}
Let $l_1=(h e^{i \varphi}, M)\, ,$ $l_2=(i, h e^{i \varphi}) \, ;$ let 
$\Delta_1$ be a triangle with vertices at the points $0,M, h e^{i \varphi}\, ,$
let $\Delta_2$ be a triangle with vertices at the points 
$0, h e^{i \varphi}, i\, ;$ let $u(x,y)$ be the potential function of the
quadrilateral $(Q( { \varphi},h); M,M+i,i, h e^{i \varphi} )$ and let
$M_1$ be the modulus of this quadrilateral; and finally let 
$\partial /\partial n$ stand for the differentiation in 
the direction orthogonal to the boundary and pointing into the domain. 
Green's formula gives
$$ 0=\int_{\partial Q( { \varphi},h)} 
\left(  y \frac{\partial u }{\partial n} - 
u \frac{\partial y }{\partial n} \right) ds 
$$
$$ =\int_{ l_1 }   y \frac{\partial u }{\partial n} ds -M_1 -
\int_{ l_2 }   u \frac{\partial y }{\partial n}  ds + M
$$
$$ =\int_{ l_1 }   y \frac{\partial (u-y) }{\partial n} ds 
+I(y,\Delta_1) -M_1 - 
\int_{ l_2 }   (u-y) \frac{\partial y }{\partial n}  ds -
I(y,\Delta_2) +M + O(h)
$$
$$
=  \frac{1}{2}  \, M\,  h \,  \sin \varphi- M_1 -\frac{1}{2} \, h
 \cos \varphi +M + O(h)
$$
when $h \to 0\,.$
\end{pf}

From Theorem \ref{th2.3} it follows that in the case of a rectangle in
a neighborhood of the vertex $z=0$ we have the situation depicted in
Figure 4.

\bigskip
\begin{figure}[ht]
\centerline{\psfig{figure=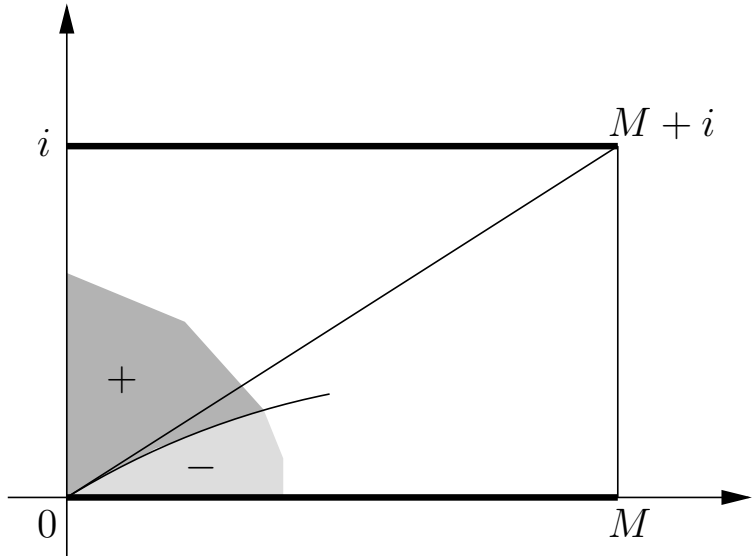,width=10cm}}
\caption{The change of modulus in a neighborhood of the vertex $z=0$
in the case of a rectangle.}
\end{figure}

\bigskip


\cc
\section{Polarization of quadrilaterals}{}

The title of this section refers to a transformation of condensers.
In the particular case of quadrilaterals it is expedient to define
this transformation directly, not via transformation of the plates as 
it was carried out in \cite{D2}, for example. Let the vertices $a,b$ of the
quadrilateral $(Q;a,b,c,d)$ be located on a ray emanating from the origin
and the vertices $c,d$ on the ray obtained from the first one by reflection
in the real axis. We also allow for the limiting case when the angle between 
the rays is zero; 
then both lines are parallel to the $x$ -axis, one of them
in the upper, the other one in the lower half plane (Figure 5).

\bigskip
\begin{figure}[ht]
\centerline{\psfig{figure=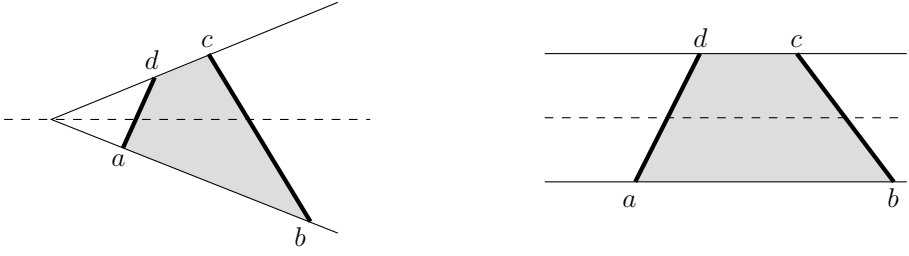,width=12cm}}
\caption{Before polarization.}
\end{figure}

We say that the result of the polarization of the quadrilateral $(Q;a,b,c,d)$
(with respect to real axis) is the quadrilateral 
$(PQ;a,\overline{c},\overline{b},d) \,.$
Figure 6 displays the resulting polarized quadrilaterals in the both cases.

\bigskip
\begin{figure}[ht]
\centerline{\psfig{figure=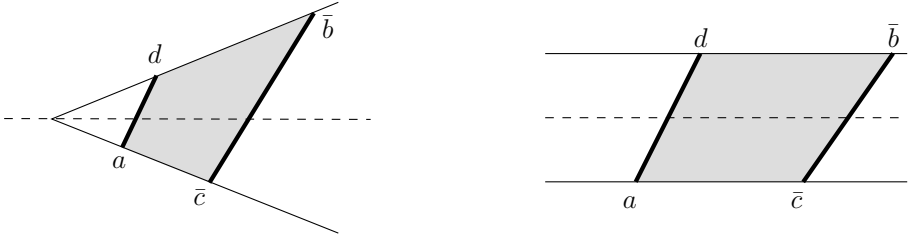,width=12cm}}
\caption{After polarization.}
\end{figure}

\begin{theo}{}\label{th3.1}
If the points $a,b,c,d$ are as above, and if 
$Re a \le Re d\, ,$ $Re c < Re b \,,$ (or $Re a < Re d\, ,$ 
$Re c \le Re b \,,$)
then the following inequalities hold
$$M(Q;a,b,c,d) > M(PQ;a,\overline{c}, \overline{b},d) \, ,$$
$$M(Q';d,c,b,a) > M(PQ';d,\overline{b}, \overline{c},a) \, .$$
\end{theo}

\begin{pf} It is enough to consider the case when
$Im a <0 \, ,$ $Im b <0 \,,$ and $Re a < Re b  \, $ (Fig. 5) .
Let $u(z)$ be the potential function of the quadrilateral $(Q;a,b,c,d) \,.$
We extend it by continuity to be equal to zero or one in the intersection
with the domain $Re a \le Re z \le Re b$, respectively, with the sector or the
horizontal strip, on the boundary of which the vertices $a,b,c,d$ are located.
Following the proof of Theorem 1.1 of \cite{D2} we introduce an auxiliary
function
$$
Pu(z)=
\begin{cases}
\min \{ u(z) , u(\overline{z}) \} , &  Im z \ge 0 \, , \\
\max \{ u(z) , u(\overline{z}) \} , &  Im z \le 0  \, , 
\end{cases}
$$
for $z \in \overline{PQ} \,.$ This function is admissible for the quadrilateral
 $(PQ;a,\overline{c},\overline{b},d) \,.$ Therefore we obtain from
(\ref{eq1.1}) and (\ref{eq1.2}) as in \cite{D2}
$$
M(Q;a,b,c,d) = I(u,Q)= I(Pu,PQ) \ge M(PQ;a,\overline{c},\overline{b},d) \,. 
$$
The equality is possible only in the case when the function $Pu(z)$ coincides
with the potential function of the quadrilateral
 $(PQ;a,\overline{c},\overline{b},d) \,.$ In this case, by virtue of uniqueness
$Pu(z) \equiv u(z), z \in Q,$ which contradicts the hypothesis of the theorem.
Thus, the first inequality of Theorem \ref{th3.1} is proved. The proof of the 
second inequality is analogous, the only difference being that the potential
function must be extended in the intersection of the set 
$\{ z : Re z \le Re d\} \cup \{z : Re z \ge Re c\}$ correspondingly with the
sector or the horizontal strip, and in the definition of the function $Pu(z)$
the places of the $\max$ and $\min \,$ must be exchanged. 
The theorem is proved.
\end{pf}

\medskip

\bigskip


\cc
\section{Continuous symmetrization}{}\label{sec4}

The symmetrization of quadrilaterals was introduced by Polya and Szeg\"o 
\cite[p. 188]{PS}. They also consider continuous symmetrization of
ring condensers with the additional property that the domains bounded
by the level curves of the potential functions are convex 
\cite[p. 200]{PS}.
Combining these notions we are led to continuous symmetrization of
quadrilaterals. Let $\alpha, \beta,\gamma, \delta, \lambda$
be real numbers with 
$\alpha< \beta, \gamma> \delta, 0\le  \lambda \le 1 \,.$
We define
$$ \alpha(\lambda)= \alpha- \lambda ( \alpha+\beta)/2, \quad
 \beta(\lambda)= \beta- \lambda ( \alpha+\beta)/2,  $$
$$ \gamma(\lambda)= \gamma- \lambda ( \gamma+\delta)/2, \quad
 \delta(\lambda)= \delta- \lambda ( \gamma+\delta)/2,  $$
$$ {\mathbb Q}(\lambda) = (Q(\lambda); 1+i \alpha(\lambda), 1+i \beta(\lambda),
i \gamma(\lambda),i \delta(\lambda)) .$$
When the parameter $\lambda$ varies from $0$ to $1$ the quadrilateral
${\mathbb Q}(0)$ continuously transforms onto the trapezoid 
${\mathbb Q}(1)$ which is the
result of the Steiner symmetrization of ${\mathbb Q}(0)$ with respect to the
real axis.

The following assertion is due to Polya and Szeg\"o.

\begin{theo}{}\label{th4.1}
The function $M{\mathbb Q}(\lambda)$ monotoneously decreases on 
the segment $[0,1] \,.$
\end{theo}

\begin{pf} It is enough to prove that 
$M{\mathbb Q}(0) \ge M{\mathbb Q}(\lambda)$ for every $\lambda \in [0,1] \,.$
Following \cite[p. 188]{PS} we translate the domain ${ Q}(0)$ 
up along the imaginary axis such that the image of the domain under 
the translation, $\widetilde{Q}(0) \,,$
does not intersect with ${ Q}(0)\,.$ 
Let $u(z)$ be the potential function of
the quadrilateral ${\mathbb Q}(0),$ and $\widetilde{u}(z)$ 
be the potential function
corresponding to the quadrilateral $\widetilde{ \mathbb Q}(0)\, .$ 
Let the function $v(z) = u(z)$ for 
$z \in Q(0), v(z) = 1- \widetilde{u}(z)$ for 
$z \in \widetilde{Q}(0)$ and let $v(z)$ be defined by 
continuity, equal to zero or one
in the strip $S= \{ z: 0 \le Re z \le 1\} \,,$ respectively. 
Finally let $v_{\lambda}(z)$
be the result of the continuous symmetrization of $v(z)$ \cite[p. 203]{PS}.
We have
$$
2M {\mathbb Q}(0) = I(u, Q(0)) + 
I(\widetilde{u}, \widetilde{Q}(0)) =I(v,S) \,.
$$
According to Polya and Szeg\"o
\begin{equation} \label{eq4.2}
I(v,S)\ge I(v_{\lambda},S) .
\end{equation}
It remains to observe that the function $v_{\lambda}$ is admissible for
two nonintersecting quadrilaterals, which are obtained from 
${\mathbb Q}(\lambda)$
by translation and therefore have moduli equal to its modulus. The proof
is complete.
\end{pf}

\medskip

\begin{rem} \label{rem1}
 For the proof of (\ref{eq4.2}) it is essential that every
line $Re z = x, 0 < x < 1,$ intersects every level curve of the function
$u(z)$ only in one point. This property of the level curve for polygonal
quadrilaterals is easily established by use of standard methods of the
theory of functions \cite[p. 221]{N}. 
\end{rem}

\medskip

\begin{rem} \label{rem2}
Polarization of a trapezoid with two opposite sides of equal length
yields a parallelogram. Therefore the superposition of continuous
symmetrization and polarization gives a continuous transformation of
a quadrilateral of the form ${\mathbb Q}(0)$ onto a parallelogram having an
area equal to the area of  $Q(0)\,.$ If this parallelogram is again
symmetrized we obtain a rectangle. At each step the modulus 
of the quadrilateral does not increase. The fact that the modulus
$M {\mathbb Q} (0)$ is greater than or equal to the modulus of the
corresponding rectangle with the same area coincides with 
the classical isoperimetric inequality. 
\end{rem}

\bigskip
\begin{figure}[ht]
\centerline{\psfig{figure=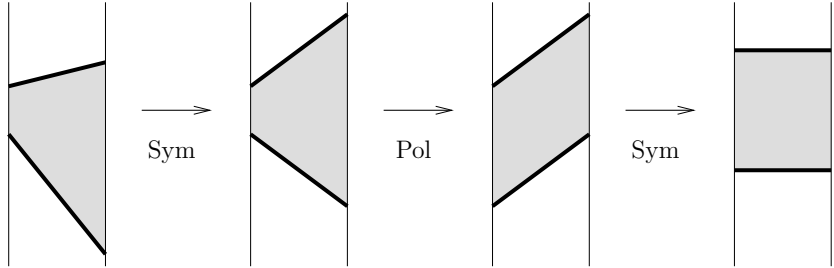,width=11cm}}
\caption{At each step the modulus 
of the quadrilateral does not increase.}
\end{figure}

In 1998 E. Reich made an interesting observation that the Polya -Szeg\"o
method yields a convex function $M{\mathbb Q}(\lambda)$ in the case when $Q(0)$
is a parallelogram \cite{R}. In fact, even a more general assertion holds,
which follows from the linear averaging transformation of Marcus \cite{M}.
\bigskip


\cc
\section{Averaging transformation and rotation of sides}{}\label{sec5}

Let $\{{\mathbb Q}_k\}_{k=1}^n$ be a collection of quadrilaterals of the form
$$
{\mathbb Q}_k = (Q_k; 1+ i \alpha_k, 1+ i \beta_k, i \gamma_k, i \delta_k)
$$
with $\alpha_k < \beta_k, \gamma_k > \delta_k, k=1,...,n,$ and let
$A=\{a_k\}, k=1,...,n,$ be a set of positive numbers with 
$\sum_{k=1}^n a_k =1 \, .$ The result of linear averaging transformations
of the family of quadrilaterals  $\{ {\mathbb Q}_k\}_{k=1}^n,$ is, 
by definition, the quadrilateral
$$
{\cal L}_A(\{ {\mathbb Q}_k\}_{k=1}^n) = (Q; 1+ i \sum_{k=1}^n a_k \alpha_k,
1+ i \sum_{k=1}^n a_k \beta_k, i \sum_{k=1}^n a_k \gamma_k,
 i \sum_{k=1}^n a_k \delta_k)\,  . 
$$

\begin{theo}{}\label{th5.1}
The following inequality holds 
$$
M {\cal L}_A(\{ {\mathbb Q}_k \}_{k=1}^n ) \le 
\sum_{k=1}^n a_k M {\mathbb Q}_k \, .
$$
\end{theo}

\begin{pf}
We may assume that $\alpha_k \ge 0\,,$  $\delta_k \ge 0$ for 
$k=1,...,n\,.$ Denote by $u_k, k=1,...,n,$ the potential function
corresponding to the quadrilateral ${\mathbb Q}_k$ extended by continuity
to the half strip
$\{z: 0 \le Re z \le 1, Im z \ge 0\}$ with values zero and one. Let
 $$
u^* = {\cal L}_A(\{ u_k\}_{k=1}^n) \,  , 
$$
be the result of the linear averaging transformation of Marcus of the family 
of functions $\{ u_k\}_{ k=1}^n \,.$ By \cite[Theorem 1.1]{M}
$$
I(u^*,Q) \le \sum_{k=1}^n a_k I(u_k,Q_k) = \sum_{k=1}^n a_k M {\mathbb Q}_k \,.
$$
It remains to observe that the function $u^*$ is admissible for the
quadrilateral ${\cal L}_A(\{ {\mathbb Q}_k\}_{k=1}^n) \, .$ 
The proof is complete.
\end{pf}

\begin{cor} \label{cor5.2}
The function $M(y) = M(Q(y); 1+i\alpha, 1+i y, i \gamma, i \delta), 
\gamma >\delta, $ is a decreasing and convex function on the interval
$(\alpha, \infty) \,.$
\end{cor}

\begin{theo}{}\label{th5.3}
For real numbers $\gamma$ and $\beta, 0 < \gamma \le \beta,$ the function
$$
q_1(y) =M(Q_1(y); 1+iy, 1+i \beta, i \gamma, -iy)
$$
has the following properties:

1) $q_1(y)$ decreases on $(-\gamma, 0)\,,$

2) for every  $y \in (-\gamma, 0)\,$ we have $q_1(y) > q_1(-y)  \, ,$

3) $q_1(y)$ is convex on  $(-\gamma, \beta)\, .$
\end{theo}

\begin{pf}
The property 3) follows from Theorem \ref{th5.1}. The property 2) is
a consequence of Theorem \ref{th3.1}. Finally, the property 1) is a
corollary to 2) and 3).
\end{pf}

Note that when $y$ increases from $-\gamma$ to $\beta$ the side 
$(-iy, 1+iy)$ of the quadrilateral $(Q_1(y); 1+iy,1+i \beta,i \gamma, -i y)$
turns around the point $z= 1/2$ counterclockwise.

\begin{theo}{}\label{th5.4}
For real numbers $\varphi, r_1, r_2$ with $0 < \varphi < \pi/2,
1/\cos \varphi <r_1 \le r_2, $ the function
$$
q_2(y) =M\,{{\mathbb Q}}(r)\, , \quad {{\mathbb Q}}(r) = 
(Q_2(r); r e^{-i \varphi}, 
r_2 e^{-i \varphi}, r_1 e^{i \varphi}, 
(2 \cos \varphi -1/r)^{-1} e^{i \varphi}   )
$$
has the following properties:

1) $q_2(y)$ decreases on $( (2 \cos \varphi -1/r_1)^{-1}, 1/ \cos \varphi)\, ,$

2) for every  $ r \in ((2 \cos \varphi -1/r)^{-1}, 
1/\cos \varphi)\, ,$ we have $q_2(r) > q_2((2 \cos \varphi -1/r)^{-1})  \, ,$

3) the function $q_2(1/p)$ is convex on  $(1/r_2, 2 \cos \varphi -1/r_1) \, .$ 
\end{theo}

\begin{pf}
Let $w= F(z)$ be a regular branch of the function 
$w= \frac{i}{2 \varphi } \ln z + \frac{1}{2} \, ,$ 
which maps the angle $|\arg z| <\varphi$ onto the strip
$0 < Re w < 1 \, .$ It is clear what is understood by the image of the
quadrilateral ${{\mathbb Q}}(r)$ under the mapping $w= F(z).$ Applying to 
the potential functions of the (nonpolygonal) quadrilaterals 
$ F {{\mathbb Q}}(r')\, ,$ $ F {{\mathbb Q}}(r'')$ 
the linear averaging transformation of
Marcus \cite{M}, in the same way as in the proof of Theorem \ref{th5.1}
we obtain the inequality
$$
M F {{\mathbb Q}}(((1/r' + 1/r'')/2)^{-1}) \le 
\frac{1}{2}(  M F {{\mathbb Q}}( r') + M F {{\mathbb Q}}(r'')) \, .
$$
This inequality reduces to property 3). Property 2) follows from Theorem
\ref{th3.1} and property 1) is a consequence of 2) and 3). The proof is
complete. 
\end{pf}

This result characterizes the change of the modulus of the quadrilateral
${{\mathbb Q}}(r)$ under the rotation of one of its sides around 
the point $z=1 \,.$


\bigskip

\cc
\section{Examples and open problems }{}\label{sec6}

The properties of the moduli of quadrilaterals studied in Sections 1-5
enable us to answer many questions about the behavior of these
moduli under geometric transformations. Complementing our earlier
examples we now make a few additional remarks.

\begin{nonsec} {Question} \label{Q3} \end{nonsec} 
Let $a,b \in C$ with $\Im a >0, \Im b >0$ and assume that $ (a,b,0,1)$
determines the vertices of a quadrilateral and
$ \arg b \in (\pi/2,\pi),$ $ \arg (a-1) \in (0,\pi/2).$
Is it true that
$$ M(Q;a,b,0,1) \le M(Q';1+i|a-1|,i|b|,0,1) ?$$

Properties 1 and 2 in Section 1 give an affirmative answer to
this question.

\begin{nonsec} {Question} \label{Q1} \end{nonsec} 
 Let $h,k >0, t \in R$ and
$$g(t) = M(Q;1+i(t+2k), ih, -ih, 1+it) .$$
Is it true that  $g$ attains its largest value when $t=-k?$


An affimative answer to this question follows from the Steiner
symmetrization (see Section 4).

\medskip

On the other hand, some problems for the moduli of quadrilaterals,
involving length and area remain open. We formulate two of them.

\begin{nonsec} {Open problem} \label{Q2} \end{nonsec} 
 Let $a,b \in C$ with $\Im a >0, \Im b >0$ and assume that $ (a,b,0,1)$
determines the vertices of a quadrilateral and $\beta\equiv \arg b >
\alpha \equiv \arg a$  Then the area of the quadrilateral $(Q;a,b,0,1)$
is $\mathrm{area}(a,b) \equiv |a|(\sin \alpha +|b|\sin(\beta-\alpha))/2.$
Let $2h=|b|, 2k =|a-1|.$

(a) Define $t$ by the condition that
$$(Q;a,b,0,1), \textrm{ and }(Q';t+ik, ih,-ih, t-ik)$$
have equal areas, i.e. that
$$ \mathrm{area}(a,b) =(|b|+|a-1|)t/2 $$
Is it true that modulus of $(Q;a,b,0,1)$ is smaller than that of 
$(Q';t+ik, ih,-ih, t-ik)$?

(b) We could also consider various other choices for $t$. For instance,
we could define $t$ as the distance of the segments $[0,b]$ and 
$[1,a]$ and ask the same question as in (a).

\begin{nonsec} {Open problem} \label{Q5} \end{nonsec} 
Let $A= 1+ r \exp(i \alpha)$, $B= s \exp (i \beta)$, with  $\arg B > \arg A$
and  $\alpha, \beta \in (0,\pi)$.
Let $$c_1=  M(Q;A,B,0,1), c_2= M(Q_1;A, A-1, 0,1), c_3= 
M(Q_3;A,B, 0, A-B)\, .$$
 
Is it true that
$$ \min \{c_2,c_3\}\le c_1 \le \max \{ c_2, c_3 \} \, ?$$
               
\bigskip
For information about the numerical computation of the modulus of a
polygonal quadrilateral see \cite{RaV}.


\bigskip

{\bf Acknowledgements.} The authors' research was supported by the 
FEBRAS grant no. 05-III-A-01-039, 
The Finnish Academy of Sciences, and the Academy of Finland. The paper
was written during the visit of the first author to the University of 
Helsinki in March 2005. 

\bigskip

 \noindent
{\bf V.N. Dubinin }\\
Institute of Applied Mathematics\\
Far East Branch Russian Academy of Sciences\\
Vladivostok\\
RUSSIA\\
Email: {\tt dubinin@iam.dvo.ru}\\

 \medskip

\noindent
{\bf M. Vuorinen}\\
Department of Mathematics\\
FIN-20014 University of Turku \\
FINLAND\\
E-mail: {\tt vuorinen@utu.fi}\\

\end{document}